\input ./style/arxiv-vmsta.cfg
\documentclass[numbers,compress,v1.0.1]{vmsta}
\usepackage{vtexbibtags}

\volume{7}% Updated by VTEXPTS2LaTeX.exe, 25.03.2020 14:30
\issue{1}% Updated by VTEXPTS2LaTeX.exe, 25.03.2020 14:30
\pubyear{2020}
\firstpage{43}% Updated by VTEXPTS2LaTeX.exe, 25.03.2020 14:30
\lastpage{60}% Updated by VTEXPTS2LaTeX.exe, 25.03.2020 14:30
\aid{VMSTA148}% Updated by VTEXPTS2LaTeX.exe, 11.02.2020 10:30
\doi{10.15559/20-VMSTA148}% Updated by VTEXPTS2LaTeX.exe, 11.02.2020 10:30
\articletype{research-article}

\DeclareFontFamily{U}{mathx}{\hyphenchar \font 45}
\DeclareFontShape{U}{mathx}{m}{n}{
<5> <6> <7> <8> <9> <10>
<10.95> <12> <14.4> <17.28> <20.74> <24.88>
mathx10
}{}
\DeclareSymbolFont{mathx}{U}{mathx}{m}{n}
\DeclareFontSubstitution{U}{mathx}{m}{n}
\DeclareMathAccent{\widecheck}{0}{mathx}{"71}

\startlocaldefs
\newtheorem{thm}{Theorem}[section]
\newtheorem{lem}{Lemma}[section]
\newtheorem{cor}{Corollary}[section]
\newtheorem{prop}{Proposition}[section]
\newtheorem*{Claim}{Claim}

\theoremstyle{definition}
\newtheorem{symb}{Notation}[section]
\newtheorem{symbs}{Notations}[section]
\newtheorem{rem}{Remark}[section]
\newtheorem{ex}{Example}[section]

\newcolumntype{d}[1]{D{.}{.}{#1}}

\allowdisplaybreaks
\ifdefined\HCode % tex4ht
\def\index#1{}
\def\texorxml#1#2{#2}
\else % LaTeX
\def\texorxml#1#2{#1}
\fi
\endlocaldefs

\begin{document}

\begin{frontmatter}
\pretitle{Research Article}

\title{A characterization of equivalent martingale measures in a renewal risk model with
applications to premium calculation principles}

\author{\inits{N.D}\fnms{Nikolaos D.}~\snm{Macheras}\ead[label=e1]{macheras@unipi.gr}}
\author{\inits{S.M.}\fnms{Spyridon M.}~\snm{Tzaninis}\thanksref{cor1}\ead[label=e2]{stzaninis@unipi.gr}}%
\thankstext[type=corresp,id=cor1]{Corresponding author.}
\address{Department of Statistics and Insurance Science, \institution{University of Piraeus},
80 Karaoli and Dimitriou \xch{Street,}{street} 185 34 Piraeus, \cny{Greece}}

%\author[]{\inits{}\fnms{}~\snm{}\thanksref{f1}\ead[label=e1]{}}
%\author[]{\inits{}\fnms{}~\snm{}\thanksref{f1}\thanksref{cor1}\ead[label=e2]{}}
%\thankstext[type=corresp,id=cor1]{Corresponding author.}
%\address[]{\institution{}, ..., \cny{}}
%\address[]{\institution{}, ..., \cny{}}

%\thankstext[id=f1]{}
\markboth{N.D. Macheras, S.M. Tzaninis}{A characterization of equivalent martingale measures for CRPs}
%\dedicated{}

%\markboth{Authors}{Title}
%\markboth{}{}
\begin{abstract}
Generalizing earlier work of Delbaen and Haezendonck for given com\-pound
renewal process $S$ under a probability measure $P$ we characterize all
probability measures $Q$ on the domain of $P$ such that $Q$ and $P$ are
progressively equivalent and $S$ remains a compound renewal process under
$Q$. As a consequence, we prove that any compound renewal process can be
converted into a compound Poisson process through a change of measures and
we show how this approach is related to premium calculation principles.
\end{abstract}



\received{\sday{19} \smonth{4} \syear{2019}}% Updated by VTEXPTS2LaTeX.exe, 11.02.2020 10:30
\revised{\sday{3} \smonth{2} \syear{2020}}% Updated by VTEXPTS2LaTeX.exe, 11.02.2020 10:30
\accepted{\sday{3} \smonth{2} \syear{2020}}% Updated by VTEXPTS2LaTeX.exe, 11.02.2020 10:30
\publishedonline{\sday{20} \smonth{2} \syear{2020}}
\end{frontmatter}

%s1 #&#
\section{Introduction}\label{intro}

A basic method in mathematical finance is to replace the original probability
measure\index{probability measure} with an equivalent martingale measure,\index{equivalent martingale measures} sometimes called a
risk-neutral measure. This measure is used for pricing and hedging given
contingent claims (e.g., options, futures, etc.). In contrast to the
situation of the classical Black--Scholes option pricing formula, where the
equivalent martingale measure\index{equivalent martingale measures} is unique, in actuarial mathematics that is
certainly not the case.

The above fact was pointed out by Delbaen and Haezendonck in their pioneering
paper \cite{dh}, as the authors ``tried to create a mathematical framework to
deal with finance related to risk processes'' in the frame of classical Risk
Theory. Thus, they were confronted with the problem of characterizing all
equivalent martingale measures $Q$\index{equivalent martingale measures} such that a compound Poisson process\index{compound Poisson ! process} under
an original measure $P$ remains a compound Poisson\index{compound Poisson} one under $Q$. They
%answered to the positive
solved positively the previous problem in \cite{dh}, and applied their
results to the theory of premium calculation principles\index{premium calculation principles} (see also Embrechts
\cite{em} for an overview). The method provided by \cite{dh} has been
successfully applied to many areas of insurance mathematics such as pricing
(re-)insurance contracts (Holtan \cite{ho}, Haslip and Kaishev \cite{haka}),
simulation of ruin probabilities (Boogaert and De Waegenaere \cite{bw}), risk
capital allocation (Yu et al. \cite{yu}), pricing CAT derivatives (Geman and
Yor \cite{gy}, Embrechts and Meister \cite{emme}), and has been generalized to
the case of %the
mixed Poisson processes (see Meister \cite{me}).

However, there is one vital point about the (compound) Poisson processes
which is their greatest weakness as far as practical applications are
considered, and this is the fact that the variance is a linear function of
time $t$. The latter, together with the fact that in some interesting
real-life cases the interarrival times process associated with a counting
process remain independent but the exponential interarrival time distribution
does not fit well into the observed data (cf. e.g. Chen et al. \cite{chen}
and Wang et al. \cite{wang}), implies that the induced counting process is a
renewal but not a Poisson one. This raises the question, whether the
characterization of Delbaen and Haezendonck can be extended to the more
general compound renewal risk\index{compound renewal risk} model (also known as the Sparre--Andersen
model), and it is precisely this problem the paper deals with. In particular,
if the process $S$ is under the probability measure\index{probability measure} $P$ a compound renewal\index{compound renewal}
one, it would be interesting to characterize all probability measures\index{probability measure} $Q$
being equivalent to $P$ and converting $S$ into a compound Poisson process\index{compound Poisson ! process}
under $Q$.

In Section~\ref{CRPPEM}, we prove the one direction of the desired
characterization, see\break Proposition~\ref{thm1}, which provides characterization
and explicit calculation of Radon--Nikod\'{y}m derivatives $dQ/dP$ for
well-known cases in insurance mathematics, see Examples~\ref{ex1b} and~\ref{ex3b}.
Since the increments of a renewal process are not, in general, independent
and stationary we cannot use arguments similar to those used in the main
proof of \xch{\cite[Proposition 2.2]{dh}}{\cite{dh}, Proposition 2.2}. In an effort to overcome this obstacle
we inserted Lemma~\ref{lem1}, which holds true for any (compound) counting
process, and on which the proof of Proposition~\ref{thm1} relies heavily.

In Section~\ref{char}, the inverse direction is proven in
Proposition~\ref{thm2}, where a canonical change of measures technique is
provided, which seems to simplify the well-known one involving the
markovization of a (compound) renewal process, see Remark~\ref{gen}. The
desired characterization is given in Theorem~\ref{thm!}, which completes and
simplifies the proof of the main result of \cite{dh}. As a consequence of
Theorem~\ref{thm!}, it is proven in Corollary~\ref{cor4b} that any compound
renewal process\index{compound renewal ! process} can be converted into a compound Poisson\index{compound Poisson} one through a change
of measures, by choosing the ``correct'' Radon--Nikod\'{y}m derivative. The
main result of \xch{\cite[Proposition~2.2]{dh}}{\cite{dh}, Proposition 2.2,}
follows as a special instance of Theorem~\ref{thm!}, see Remark~\ref{cla}
(a).

In Section~\ref{CRPM}, we apply our results to the financial pricing of
insurance in a compound renewal risk\index{compound renewal risk} model. We first prove that given a
compound renewal process\index{compound renewal ! process} $S$ under $P$, the process $Z(P):=\{Z_{t}\}_{t\in
\mathbb{R}_{+}}$ with $Z_{t}:=S_{t}-t\cdot p(P)$ for any $t\geq 0$, where
$p(P)$ is the premium density,\index{premium density} is a martingale under $P$ if and only if $S$
is a compound Poisson process\index{compound Poisson ! process} under $P$, see Proposition~\ref{thm4}, showing
in this way that a martingale approach to premium calculation principles
leads\index{premium calculation principles} in the case of compound renewal processes\index{compound renewal ! process} immediately to compound
Poisson\index{compound Poisson} ones. A consequence of Theorem~\ref{thm!} and Proposition~\ref{thm4}
is a characterization of all progressively equivalent martingale measures\index{equivalent martingale measures} $Q$
converting a compound renewal process\index{compound renewal ! process} $S$ into a compound Poisson\index{compound Poisson} one, see
Proposition~\ref{mar}. Using the latter result, we find out \emph{canonical}
price processes satisfying the condition of no free lunch with vanishing
risk, see Theorem~\ref{nfl}, connecting in this way our results with this
basic notion of mathematical finance. Finally, we present some applications
of Corollary~\ref{cor4b} and Theorem~\ref{nfl} to the computation of some
premium calculation principles,\index{premium calculation principles} see Examples~\ref{expcp0} to~\ref{expcp7}.

%s2 #&#
\section{Compound renewal processes\index{compound renewal ! process} and progressively equivalent measures}\label{CRPPEM}

\emph{Throughout this paper, unless stated otherwise, $(\varOmega ,\varSigma
,P)$ is a fixed but arbitrary probability space and $\varUpsilon :=(0,\infty
)$.} The symbols $\mathcal{L}^{1}(P)$ and $\mathcal{L}^{2}(P)$ stand for the
families of all real-valued $P$-integrable and $P$-square integrable
functions on $\varOmega $, respectively. Functions that are $P$-a.s. equal
are not identified. We denote by $\sigma (\mathcal{G})$ the $\sigma $-algebra
generated by a family $\mathcal{G}$ of subsets of $\varOmega $. Given a
topology $\mathfrak{T}$ on $\varOmega $ we write ${\mathfrak{B}}(\varOmega )$
for its \textbf{Borel $\sigma $-algebra} on $\varOmega $, i.e. the $\sigma
$-algebra generated by $\mathfrak{T}$. Our measure theoretic terminology is
standard and generally follows \cite{Co}. For the definitions of real-valued
random variables and random variables we refer to \xch{\cite[p.~308]{Co}}{\cite{Co}, p.~308}. We apply
the notation $P_{X}:=P_{X}(\theta ):=\mathbf{{K}}(\theta )$ to mean that $X$
is distributed according to the law $\mathbf{{K}}(\theta )$, where $\theta
\in D\subseteq \mathbb{R}^{d}$ ($d\in \mathbb{N}$) is the parameter of the
distribution. We denote again by $\mathbf{K}(\theta )$ the distribution
function induced by the probability distribution $\mathbf{K}(\theta )$.
Notation $\mathbf{Ga}(a,b)$, where $a,b\in (0,\infty )$, stands for the law
of gamma distribution (cf. e.g. \xch{\cite[p.~180]{Sc}}{\cite{Sc}, p.~180}). In particular,
$\mathbf{Ga}(a,1)=\mathbf{Exp}(a)$ stands for the law of exponential
distribution. For two real-valued random variables $X$ and $Y$ we write $X=Y$
$P$-a.s. if $\{X\neq Y\}$ is a $P$-null set. If $A\subseteq \varOmega $, then
$A^{c}:=\varOmega \setminus A$, while $\chi _{A}$ denotes the indicator (or
characteristic) function of the set $A$. For a map $f:D\rightarrow E$ and for
a nonempty set $A\subseteq D$ we denote by $f\upharpoonright A$ the
restriction of $f$ to $A$. For the unexplained terminology of Probability and
Risk Theory we refer to \cite{Sc}.

A sequence $W:=\{W_{n}\}_{n\in \mathbb{N}}$ of positive real-valued random
variables on $\varOmega $ is called a \textbf{(claim) interarrival process}
(cf. e.g. \xch{\cite[p.~7]{Sc}}{\cite{Sc}, p.~7}). The \textbf{(claim)
arrival process} $T:=\{T_{n}\}_{n\in \mathbb{N}_{0}}$ induced by $W$ is
defined by means of $T_{0}:=0$ and $T_{n}:=\sum_{k=1}^{n} W_{k}$ for any
$n\in \mathbb{N}$ (cf. e.g. \xch{\cite[p.~7]{Sc}}{\cite{Sc}, p.~7}). A
\textbf{counting} (or \textbf{claim number}) \textbf{process}
$N:=\{N_{t}\}_{t\in \mathbb{R}_{+}}$ is defined by means of
$N_{t}:=\sum_{n=1}^{\infty } \chi _{\{T_{n}\leq t\}}$ for any $t\geq 0$ (cf.
e.g. \xch{\cite[Theorem 2.1.1]{Sc}}{\cite{Sc}, Theorem 2.1.1}). In
particular, if $W$ is $P$-i.i.d. with common distribution $\mathbf{K}(\theta
):\mathfrak{B}(\varUpsilon )\rightarrow [0,1]$ ($ \theta \in D\subseteq
\mathbb{R}^{d}$), the counting process $N$ is a $P$-\textbf{renewal process
with parameter $\theta \in D\subseteq \mathbb{R}^{d}$ and interarrival time
distribution $\mathbf{K}(\theta )$} (written $P$-RP$(\mathbf{K}(\theta ))$
for short). If $\theta >0$ and $\mathbf{K}(\theta )=\mathbf{Exp}(\theta )$
then a $P$-RP$(\mathbf{K}(\theta ))$ becomes a $P$-Poisson process with
parameter $\theta $ (cf. e.g. \xch{\cite[p.~23 for the
definition]{Sc}}{\cite{Sc}, p.~23 for the definition}). Note that if $N$ is a
$P$-RP$(\mathbf{K}(\theta ))$ then $\mathbb{E}_{P}[N^{m}_{t}]<\infty $ for
any $t\geq 0$ and $m\in \mathbb{N}$ (cf. e.g. \xch{\cite[Proposition 4,
p.~101]{se}}{\cite{se}, Proposition 4, p.~101}); hence according to
\xch{\cite[Corollary~2.1.5]{Sc}}{\cite{Sc}, Corollary~2.1.5}, it has zero probability of explosion, i.e.
$P(\{\sup_{n\in \mathbb{N}} T_{n}<\infty \})=0$. Furthermore, if
$X:=\{X_{n}\}_{n\in \mathbb{N}}$ is another sequence of
$P$-\xch{i.i.d.}{i.i.d} positive real-valued random variables on $\varOmega
$, called \textbf{claim size process} (cf. e.g. \xch{\cite[p.~103]{Sc}}{\cite{Sc}, p.~103}), which is
independent of $N$, define the \textbf{aggregate claims process}
$S:=\{S_{t}\}_{t\in \mathbb{R}_{+}}$ by means of $S_{t}:=\sum_{n=1}^{N_{t}}
X_{n}$ for any $t\geq 0$ (cf. e.g. \xch{\cite[p.~103]{Sc}}{\cite{Sc}, p.~103}). In particular, if $N$
is a $P$-RP$(\mathbf{K}(\theta ))$, the aggregate claims process is a
$P$-\textbf{{compound renewal process}} ($P$-CRP for short) \textbf{{with
parameters $\mathbf{K}(\theta )$ and $P_{X_{1}}$}}. In the special case where
$N$ is a $P$-Poisson process with parameter $\theta $, the aggregate claims
process $S$ is called a $P$-\textbf{{compound Poisson process}} ($P$-CPP for
short) \textbf{{with parameters $\theta $ and $P_{X_{1}}$}}.

\emph{Henceforth, unless stated otherwise, $S:=\{S_{t}\}_{t\in
\mathbb{R}_{+}}$ is a $P$-CRP with parameters $\mathbf{K}(\theta )$ and
$P_{X_{1}}$, $\mathcal{F}^{W}:=\{\mathcal{F}^{W}_{n}\}_{n\in \mathbb{N}}$,
$\mathcal{F}^{X}:=\{\mathcal{F}^{X}_{n}\}_{n\in \mathbb{N}}$ and
$\mathcal{F}:=\{\mathcal{F}_{t}\}_{t\in \mathbb{R}_{+}}$ are the natural
filtrations of $W$, $X$ and $S$, respectively.}

The next lemma is a general and helpful result, as it provides a clear
understanding of the structure of $\mathcal{F}$, and it is essential for the
proofs of our main results. Lemma 2.1 is a part of \xch{\cite[Lemma III.1.29]{js}}{\cite{js}, Lemma III.1.29},
but we write it with its proof in a form suitable for our results.

%l2.1 #&#
\begin{lem}\label{lem1}
For every $t\geq 0$ and $n\in \mathbb{N}_{0}$ the following
\begin{equation*}
\mathcal{F}_{t}\cap \{N_{t}=n\}= \sigma \bigl(
\mathcal{F}^{W}_{n}\cup \mathcal{F}^{X}_{n}
\bigr)\cap \{N_{t}=n\}
\end{equation*}
holds true.
\end{lem}

\begin{proof}
Fix an arbitrary $t\geq 0$ and $n\in \mathbb{N}_{0}$.

Clearly, for $n=0$ we get
$\mathcal{F}^{X}_{0}=\mathcal{F}^{W}_{0}=\{\emptyset ,\varOmega \}$ and
$\mathcal{F}_{t}\cap \{N_{t}=0\}=\{\emptyset ,\varOmega \}\cap \{N_{t}=0 \}$;
hence $\mathcal{F}_{t}\cap \{N_{t}=0\}=\sigma (\mathcal{F}^{W}_{0}\cup
\mathcal{F}^{X}_{0})\cap \{N_{t}=0\}$.

\textbf{(a)} Inclusion $\sigma (\mathcal{F}^{W}_{n}\cup
\mathcal{F}^{X}_{n})\cap \{N_{t}=n\} \subseteq \mathcal{F}_{t}\cap
\{N_{t}=n\}$ holds true.

To show (a), fix an arbitrary $k\in \{1,\ldots ,n\}$. Note that $S$ is
progressively measurable with respect to $\mathcal{F}$ (cf. e.g.
\xch{\cite[p.~4 for the definition]{ks}}{\cite{ks}, p.~4 for the
definition}), since $S_{t}$ is $\mathcal{F}_{t}$-measurable and has right
continuous paths (cf. e.g. \xch{\cite[Proposition~1.13]{ks}}{\cite{ks},
Proposition~1.13}). The latter, together with the fact that $T_{k}$ is a
stopping time of $\mathcal{F}$, implies that $S_{T_{k}}$ is
$\mathcal{F}_{T_{k}}$-measurable, where $\mathcal{F}_{T_{k}}:=\{A\in
\varSigma : A\cap \{T_{k}\leq v\}\in \mathcal{F}_{v}\,\,\text{for any}\,\,
v\geq 0\}$ (cf. e.g. \xch{\cite[Proposition~2.18]{ks}}{\cite{ks},
Proposition~2.18}). But $T_{k-1}<T_{k}$ yields
$\mathcal{F}_{T_{k-1}}\subseteq \mathcal{F}_{T_{k}}$ (cf. e.g.
\xch{\cite[Lemma~2.15]{ks}}{\cite{ks}, Lemma 2.15}), implying that
$S_{T_{k-1}}$ is $\mathcal{F}_{T_{k}}$-measurable. Consequently, the random
variable $X_{k}=S_{T_{k}}-S_{T_{k-1}}$ is $\mathcal{F}_{T_{k}}$-measurable;
hence $\mathcal{F}^{X}_{n}\cap \{N_{t}=n\}\subseteq \mathcal{F}_{t}\cap
\{N_{t}=n \}$.

Since $W_{k}$ is $\mathcal{F}_{T_{k}}$-measurable, standard computations
yield $\mathcal{F}^{W}_{n}\cap \{N_{t}=n\}\subseteq \mathcal{F}_{t}\cap
\{N_{t}=n \}$, completing in this way the proof of (a).

\textbf{(b)} Inclusion $\mathcal{F}_{t}\cap \{N_{t}=n\}\subseteq \sigma
(\mathcal{F}^{W}_{n} \cup \mathcal{F}^{X}_{n})\cap \{N_{t}=n\}$ holds true.

To show (b), let $A\in \bigcup_{u\leq t}\sigma (S_{u})$. There exist an
index $u\in [0,t]$ and a set $B\in \mathfrak{B}(\varUpsilon )$ such that
$A=S_{u}^{-1}(B)=\bigcup_{m\in \mathbb{N}_{0}} (\{N_{u}=m\}\cap B_{m})$,
where $B_{m}:=(\sum^{m}_{j=1} X_{j})^{-1}(B)\in \mathcal{F}^{X}_{m}$ for any
$m\in \mathbb{N}_{0}$, implying
\begin{equation*}
A\cap \{N_{t}=n\}=D_{n}\cap \{N_{t}=n\}),
\end{equation*}
where $D_{n}:=  (\bigcup^{n-1}_{m=0} (\{N_{u}=m\}\cap B_{m})  )\cup (
\{T_{n}\leq u\}\cap B_{n})\in \sigma (\mathcal{F}^{W}_{n}\cup
\mathcal{F}^{X}_{n})$; hence $\bigcup_{u\leq t}\sigma (S_{u})\cap
\{N_{t}=n\}\subseteq \sigma ( \mathcal{F}^{W}_{n}\cup
\mathcal{F}^{X}_{n})\cap \{N_{t}=n\}$, implying (b).
This completes the proof of the lemma.
\end{proof}

%l2.2 #&#
\begin{lem}\label{lem2}
Let $Q$ be a probability measure\index{probability measure} on $\varSigma $.

\textup{\textbf{(a)}} If $X$ is $Q$-i.i.d., $Q_{X_{1}}\sim P_{X_{1}}$ and $h$
is a real-valued, one-to-one, $\mathfrak{B}(\varUpsilon )$-measurable
function, then there exists a $P_{X_{1}}$-a.s. unique real-valued
$\mathfrak{B}(\varUpsilon )$-measur\-able function $\gamma $ such that
\begin{enumerate}
\item[$(i)$] $\mathbb{E}_{P}  [h^{-1}\circ \gamma \circ X_{j}
    ]=1$;
\item[$(ii)$] for every $n\in \mathbb{N}_{0}$ and for all $A\in
    \mathcal{F}^{X}_{n}$ the condition
%
%e1 #&#
\begin{equation}
Q(A)=\mathbb{E}_{P} \Biggl[\chi _{A}\cdot \prod
_{j=1}^{n}\,\bigl(h^{-1}\circ
\gamma \circ X_{j}\bigr) \Biggr] \label{lem23b}
\end{equation}
holds true.
\end{enumerate}

{\textup{\textbf{(b)}}} If $W$ is $Q$-i.i.d. and $Q_{W_{1}}\sim P_{W_{1}}$,
then there exists a $P_{W_{1}}$-a.s. unique positive function $r\in
\mathcal{L}^{1}(P_{W_{1}})$ such that for every $n\in \mathbb{N}_{0}$ and for
all $D\in \mathcal{F}^{W}_{n}$ the condition
\begin{equation*}
Q(D)=\mathbb{E}_{P} \Biggl[\chi _{D}\cdot \prod
_{j=1}^{n}\,(r\circ W_{j})
\Biggr]
\end{equation*}
holds true.
\end{lem}

\begin{proof}
For \textbf{(a)}: First note that $h(\varUpsilon ):=\{h(y) : y\in \varUpsilon
\}\in \mathfrak{B}( \mathbb{R})$ (cf. e.g. \xch{\cite[Theorem~8.3.7]{Co}}{\cite{Co} Theorem 8.3.7}) and that
the function $h^{-1}$ is $\mathfrak{B}(h(\varUpsilon
))$-$\mathfrak{B}(\varUpsilon )$-measurable (cf. e.g. \xch{\cite[Proposition~8.3.5]{Co}}{\cite{Co} Proposition~8.3.5}). Since $P_{X_{1}}\sim Q_{X_{1}}$, by the Radon--Nikod\'{y}m
Theorem there exists a positive Radon--Nikod\'{y}m derivative $f\in
\mathcal{L}^{1}(P_{X_{1}})$ of $Q_{X_{1}}$ with respect to $P_{X_{1}}$. Put
$\gamma :=h\circ f$. An easy computation justifies the validity of $(i)$.

To check the validity of $(ii)$, fix an arbitrary $n\in \mathbb{N}_{0}$
and consider the family $\mathcal{C}_{n}:=  \{  \bigcap^{n}_{j=1} A_{j}
: A_{j}\in \sigma (X_{j})   \}  $. Standard computations show that any
$A\in \mathcal{C}_{n}$ satisfies condition \eqref{lem23b}. By a monotone
class\index{monotone class} argument it can be shown that \eqref{lem23b} remains valid for any
$A\in \mathcal{F}^{X}_{n}$.

Applying similar arguments as above we obtain (b).
\end{proof}

%s2.1 #&#
\begin{symbs}\label{symb1}
\textbf{(a)} Let $h$ be a function as in Lemma~\ref{lem2}. The class of all
real-valued $\mathfrak{B}(\varUpsilon )$-measurable functions $\gamma $ such
that $\mathbb{E}_{P}  [h^{-1}\circ \gamma \circ X_{1}  ]=1$ will be
denoted by $\mathcal{F}_{P,h}:=\mathcal{F}_{P,X_{1},h}$.

\textbf{(b)} Let us fix an arbitrary $\theta \in D\subseteq
\mathbb{R}^{d}$ and let $\boldsymbol{\Lambda }(\widetilde{\theta })$ be a
probability distribution on $\mathfrak{B}(\varUpsilon )$, where
$\widetilde{\theta }:=\rho (\theta )$ is a parameter depending on $\theta $
and $\rho $ is a function from $D$ into $\mathbb{R}^{k}$ ($d, k\in
\mathbb{N}$). The class of all probability measures\index{probability measure} $Q$ on $\varSigma $ being
\textbf{progressively equivalent} to $P$, i.e. $Q\upharpoonright
\mathcal{F}_{t}\sim P\upharpoonright \mathcal{F}_{t}$ for any $t\geq 0$, and
$S$ is a $Q$-CRP with parameters $\boldsymbol{\Lambda }(\widetilde{\theta })$ and
$Q_{X_{1}}$ will be denoted by $\mathcal{M}_{S,\boldsymbol{\Lambda
}(\widetilde{\theta })}$. In the special case $d=k$ and $\rho :=id_{D}$ we
write $\mathcal{M}_{S,\boldsymbol{\Lambda }(\theta
)}:=\mathcal{M}_{S,\boldsymbol{\Lambda }( \widetilde{\theta })}$, for simplicity.
\end{symbs}

\emph{From now on, unless stated otherwise, $h$ is a function as in
Lemma~\ref{lem2}, $D$, $\theta $ and $\widetilde{\theta }$ are as in
Notation~\ref{symb1} (b)}.

For the definition of a $(P,\mathcal{Z})$\textbf{-martingale}, where
$\mathcal{Z}:=\{\mathcal{Z}_{t}\}_{t\in \mathbb{R}_{+}}$ is a filtration on
$(\varOmega ,\varSigma )$ we refer to \xch{\cite[p.~25]{Sc}}{\cite{Sc}, p.~25}. A
$(P,\mathcal{Z})$-martingale $\{Z_{t}\}_{t\in \mathbb{R}_{+}}$ is
\textbf{$P$-a.s. positive} if $Z_{t}$ is $P$-a.s. positive for each $t\geq
0$. For $\mathcal{Z}=\mathcal{F}$ we write $P$-martingale instead of
$(P,\mathcal{F})$-martingale, for simplicity.

For a given aggregate claims process $S$ on $(\varOmega ,\varSigma )$, in
order to investigate the existence of progressively equivalent martingale
measures\index{equivalent martingale measures} (see Section~\ref{CRPM}), one has to be able to characterize
Radon--Nikod\'{y}m derivatives $dQ/dP$. Proposition~\ref{thm1} follows also as
a special case of \xch{\cite[Proposition 4.3 and Theorem 5.1]{ja}}{\cite{ja}, Proposition 4.3 and Theorem 5.1}, but we write it
in a form suitable for our purposes, and we present a rather elementary
proof.

%p2.1 #&#
\begin{prop}\label{thm1}
Let $Q$ be a probability measure\index{probability measure} on $\varSigma $ such that $S$ is a $Q$-CRP
with parameters $\boldsymbol{\Lambda }(\widetilde{\theta })$ and $Q_{X_{1}}$.
Then the following are equivalent:
\begin{enumerate}
\item[$(i)$] $Q\upharpoonright \mathcal{F}_{t}\sim P\upharpoonright
    \mathcal{F}_{t}$ for any $t\geq 0$;
\item[$(ii)$] $Q_{X_{1}}\sim P_{X_{1}}$ and $Q_{W_{1}}\sim P_{W_{1}}$;
\item[$(iii)$] there exists a $P_{X_{1}}$-a.s. unique function $\gamma \in
    \mathcal{F}_{P,h}$ such that
%
%e{RRM} #&#
\begin{equation}
\tag{RRM} Q(A)=\int_{A} M^{(\gamma )}_{t}(
\theta )\,dP\quad \text{for all } A\in \mathcal{F}_{t},
\label{mart}
\end{equation}
with
\begin{equation*}
M^{(\gamma )}_{t}(\theta ):= \Biggl[\prod
_{j=1}^{N_{t}}\,\bigl(h^{-1}\circ \gamma
\bigr) (X_{j})\cdot \frac{dQ_{W_{1}}}{dP_{W_{1}}}(W_{j}) \Biggr]
\cdot \frac{1-\boldsymbol{{\Lambda }}(\widetilde{\theta
})(t-T_{N_{t}})}{1-\mathbf{{K}}(\theta )(t-T_{N_{t}})},
\end{equation*}
where the family $M^{(\gamma )}(\theta ):=\{M^{(\gamma )}_{t}(\theta
)\}_{t\in \mathbb{R}_{+}}$ is a $P$-a.s. positive $P$-martingale satisfying
the condition $\mathbb{E}_{P}[M^{(\gamma )}_{t}(\theta )]=1$.
\end{enumerate}
\end{prop}

\begin{proof}
Fix an arbitrary $t\geq 0$.

For $(i)\Longrightarrow (ii)$: Statement $Q_{X_{1}}\sim P_{X_{1}}$ follows by
\xch{\cite[Lemma 2.1]{dh}}{\cite{dh}, Lemma 2.1}. To show statement $Q_{W_{1}}\sim P_{W_{1}}$, let $B\in
\mathfrak{B}(\varUpsilon )$ such that $Q_{W_{1}}(B)=0$. Since
$P_{W_{1}}(B)=\lim_{m\rightarrow \infty } P(W_{1}^{-1}(B)\cap \{T_{1} \leq
m\})$, $W_{1}^{-1}(B)\cap \{T_{1}\leq m\}\in \mathcal{F}_{m}$ and
$Q\upharpoonright \mathcal{F}_{m}\sim P\upharpoonright \mathcal{F}_{m}$ we
get $P(W_{1}^{-1}(B)\cap \{T_{1}\leq m\})=0$ for any $m\in \mathbb{N}$,
implying that $P_{W_{1}}(B)=0$. Replacing $Q$ by $P$ leads to $Q_{W_{1}}\sim
P_{W_{1}}$.

For $(ii)\Longrightarrow (iii)$: Let $A\in \mathcal{F}_{t}$ be given. By
Lemma~\ref{lem1}, for every $k\in \mathbb{N}_{0}$ there exists a set
$B_{k}\in \sigma (\mathcal{F}^{W}_{k}\cup \mathcal{F}^{X}_{k})$ such that
$A\cap \{N_{t}=k\}=B_{k}\cap \{N_{t}=k\}$. Thus, due to the fact that $N$ has
zero probability of explosion, we get
%
%e2 #&#
\begin{eqnarray}
Q(A)&=& \sum^{\infty }_{k=0} Q
\bigl(B_{k}\cap \{N_{t}=k\} \bigr)
\nonumber
\\
&=&\sum^{\infty }_{k=0} Q
\bigl(B_{k}\cap \{T_{k}\leq t\}\cap
\{W_{k+1}>t-T_{k} \}\bigr). \label{1.1}
\end{eqnarray}

Fix an arbitrary $n\in \mathbb{N}_{0}$ and put $G:=\bigcap_{j=1}^{n}
(W_{j}^{-1}(E_{j})\cap X_{j}^{-1}(F_{j}))\cap \{W_{n+1}>t-T_{n}\}$ where
$E_{j},F_{j}\in \mathfrak{B}(\varUpsilon )$ for any $j\in \{1,\ldots ,n\}$.
Then the set $G$ satisfies the condition
%
%e3 #&#
\begin{equation}
Q(G)=\int_{G} \Biggl[\prod^{n}_{j=1}
\,\bigl(h^{-1}\circ \gamma \bigr) (X_{j}) \cdot
\frac{dQ_{W_{1}}}{dP_{W_{1}}}(W_{j}) \Biggr]\cdot \frac{1-\boldsymbol{\Lambda
}(\widetilde{\theta })(t-T_{n})}{1-\mathbf{K}(\theta )(t-T_{n})} \,dP.
\label{1.2}
\end{equation}
In fact, by Lemma~\ref{lem2} and Fubini's Theorem we get
\begin{align*}
Q(G) &=\int \Biggl[\prod^{n}_{j=1}\,
\chi _{F_{j}}(x_{j})\cdot \chi _{E_{j}}(w_{j})
\cdot \bigl(h^{-1}\circ \gamma \bigr) (x_{j})\cdot
r(w_{j}) \Biggr]\cdot \frac{Q(\{W_{n+1}>t-w\})}{P(\{W_{n+1}>t-w\})}\cdot
\\
&\qquad P\bigl(\{W_{n+1}>t-w\}\bigr)\,P_{X_{1},\ldots ,X_{n}; W_{1},\ldots W_{n}}
\bigl(d(x_{1}, \ldots ,x_{n};w_{1},\ldots
,w_{n})\bigr)
\\
&=\int \chi _{G}\cdot \Biggl[\prod^{n}_{j=1}
\,\bigl(h^{-1}\circ \gamma \bigr) (X_{j}) \cdot
r(W_{j}) \Biggr]\cdot \frac{1-\boldsymbol{\Lambda }(\widetilde{\theta })(t-T_{n})}{1-\mathbf{K}(\theta )(t-T_{n})} \,dP,
\end{align*}
where $w:=\sum_{j=1}^{n} w_{j}$ and
$r(w_{j}):=\frac{dQ_{W_{1}}}{dP_{W_{1}}}(w_{j})$ for any $j\in \{1,\ldots
,n\}$; hence condition \eqref{1.2} follows. By a monotone class\index{monotone class} argument it
can be shown that \eqref{1.2} remains valid for any $C\in \sigma
(\mathcal{F}^{W}_{n}\cup \mathcal{F}^{X}_{n})\cap \{W_{n+1}>t-T_{n} \}$.

But since $B_{k}\cap \{T_{k}\leq t\}\in \sigma (\mathcal{F}^{W}_{k}\cup
\mathcal{F}^{X}_{k})$ for any $k\in \mathbb{N}_{0}$, conditions \eqref{1.1}
and \eqref{1.2} imply
\begin{equation*}
Q(A)=\sum^{\infty }_{k=0}
\mathbb{E}_{P} \Biggl[\chi _{A\cap \{N_{t}=k\}} \cdot \Biggl[\prod
^{N_{t}}_{j=1}\,\bigl(h^{-1}\circ
\gamma \bigr) (X_{j})\cdot r(W_{j}) \Biggr] \cdot
\frac{1-\boldsymbol{\Lambda }(\widetilde{\theta
})(t-T_{N_{t}})}{1-\mathbf{K}(\theta )(t-T_{N_{t}})} \Biggr].
\end{equation*}
Thus,
%
%e4 #&#
\begin{equation}
Q(A)=\mathbb{E}_{P}\bigl[\chi _{A}\cdot
M^{(\gamma )}_{t}(\theta )\bigr]\quad \text{for
all } A\in
\mathcal{F}_{t}, \label{new}
\end{equation}
implying
\begin{equation*}
\int_{A} M^{(\gamma )}_{u}(\theta )
\,dP=\int_{A} M^{(\gamma )}_{t}( \theta )
\,dP\quad \text{for all}\,\,u\in [0,t]\,\,\text{and}\,\, A\in
\mathcal{F}_{u};
\end{equation*}
hence $M^{(\gamma )}(\theta )$ is a $P$-martingale. The latter together with
condition \eqref{new} proves condition \eqref{mart}.

By \eqref{mart} for $A=\varOmega $ we obtain
\begin{equation*}
\mathbb{E}_{P}\bigl[M^{(\gamma )}_{t}(\theta )
\bigr]=\int_{\varOmega }M^{(\gamma
)}_{t}( \theta
)\,dP=Q(\varOmega )=1.
\end{equation*}
Note that $ \frac{1-\boldsymbol{\Lambda }(\widetilde{\theta
})(t-T_{N_{t}})}{1-\mathbf{K}(\theta )(t-T_{N_{t}})}$ is $P$-a.s. positive.
The latter, together with the fact that $h^{-1}\circ \gamma $ and $r$ are
$P_{X_{1}}$- and $P_{W_{1}}$-a.s. positive functions, respectively, implies
$P(\{M^{(\gamma )}_{t}(\theta )>0\})=1$.

The implication \textit{(iii)}$\Longrightarrow $\textit{(i)} is immediate.
\end{proof}

Proposition~\ref{thm1} allows us to explicitly calculate Radon--Nikod\'{y}m
derivatives for the most important insurance risk processes, as the following
two examples illustrate. In the first example we consider the case of the
Poisson process with parameter $\theta $.

%e2.1 #&#
\begin{ex}\label{ex1b}
Take $h:=\ln $, $\theta ,\widetilde{\theta }\in D:=\varUpsilon $, and let
$P\in \mathcal{M}_{S,\mathbf{Exp}(\theta )}$ and $Q\in
\mathcal{M}_{S,\mathbf{Exp}(\widetilde{\theta })}$. By Proposition~\ref{thm1}
there exists a $P_{X_{1}}$-a.s. unique function $\gamma \in
\mathcal{F}_{P,\ln }$ defined by means of $\gamma :=\ln f$, where $f$ is a
Radon--Nikod\'{y}m derivative of $Q_{X_{1}}$ with respect to $P_{X_{1}}$, such
that for all $A\in \mathcal{F}_{t}$
\begin{equation*}
Q(A)=\int_{A} e^{\sum ^{N_{t}}_{j=1}\gamma (X_{j})}\cdot \biggl(
\frac{\widetilde{\theta }}{\theta } \biggr)^{N_{t}}\cdot e^{-t\cdot (
\widetilde{\theta }-\theta )}\,dP.
\end{equation*}
\end{ex}

In our next example we consider a renewal process with gamma distributed
interarrival times.

%e2.2 #&#
\begin{ex}\label{ex3b}
Assume that $h:=\ln $, $\theta =(\xi _{1},\kappa _{1})\in D:=\varUpsilon
\times \mathbb{N}$, $\widetilde{\theta }=(\xi _{2},\kappa _{2})\in D$, and
let $P\in \mathcal{M}_{S,\mathbf{Ga}(\theta )}$ and $Q\in
\mathcal{M}_{S,\mathbf{Ga}(\widetilde{\theta })}$. By Proposition~\ref{thm1}
there exists a $P_{X_{1}}$-a.s. unique function $\gamma \in
\mathcal{F}_{P,\ln }$ such that for all $A\in \mathcal{F}_{t}$
\begin{eqnarray*}
Q(A)&=&\int_{A} e^{\sum ^{N_{t}}_{j=1}\gamma (X_{j})}\cdot \biggl(
\frac{\xi _{2}^{\kappa _{2}}\cdot \Gamma (\kappa _{1})}{\xi _{1}^{\kappa
_{1}}\cdot \Gamma (\kappa _{2})} \biggr)^{N_{t}}\cdot e^{-t\cdot (\xi
_{2} -\xi _{1})}\cdot
\frac{\sum_{i=0}^{\kappa _{2}-1}\frac{  (\xi
_{2}\cdot (t-T_{N_{t}})  )^{i}}{i!}}{\sum_{i=0}^{\kappa
_{1}-1}\frac{  (\xi _{1}\cdot (t-T_{N_{t}})  )^{i}}{i!}}
\\
&&\qquad \qquad \cdot \prod^{N_{t}}_{j=1}
W_{j}^{\kappa _{2}-\kappa
_{1}} \,dP.
\end{eqnarray*}
\end{ex}

%s3 #&#
\section{The characterization}\label{char}

We know from Proposition~\ref{thm1} that under the weak conditions
$Q_{X_{1}}\sim P_{X_{1}}$ and $Q_{W_{1}}\sim P_{W_{1}}$, the measures $P$ and
$Q$ are equivalent on each $\sigma $-algebra $\mathcal{F}_{t}$, a result that
does not, in general, hold true for $\mathcal{F}_{\infty }:=\sigma
(\bigcup_{t\in \mathbb{R}_{+}} \mathcal{F}_{t}  )$. Let us start with the
following helpful lemma.

%l3.1 #&#
\begin{lem}\label{lem3}
The following holds true
\begin{equation*}
\mathcal{F}_{\infty }=\mathcal{F}^{(W,X)}_{\infty }:=\sigma
\biggl(\bigcup_{n \in \mathbb{N}_{0}}\mathcal{F}^{W}_{n}
\cup \bigcup_{n\in \mathbb{N}_{0}} \mathcal{F}^{X}_{n}
\biggr).
\end{equation*}
\end{lem}

\begin{proof}
Inclusion $\mathcal{F}_{\infty }\subseteq \mathcal{F}^{(W,X)}_{\infty }$
follows immediately by Lemma~\ref{lem1} and the fact that $N$ has zero
probability of explosion.

To check the validity of the inverse inclusion, fix an arbitrary $n\in
\mathbb{N}_{0}$. Since $X_{n}$ is $\mathcal{F}_{T_{n}}$-measurable, we get
$X_{n}^{-1}(B)\cap \{T_{n}\leq \ell \}\in \mathcal{F}_{\infty }$ for all
$B\in \mathfrak{B}(\varUpsilon )$ and $\ell \in \mathbb{N}_{0}$; hence
$X_{n}^{-1}(B) \in \mathcal{F}_{\infty }$, implying together with the
$\mathcal{F}_{\infty }$-measurability of $T_{n}$ that
$\mathcal{F}^{(W,X)}_{\infty }\subseteq \mathcal{F}_{\infty }$.
\end{proof}

Note that the above lemma remains true, without the assumption $P\in
\mathcal{M}_{S,\mathbf{K}(\theta )}$, under the weaker assumption that $N$
has zero probability of explosion.

%r3.1 #&#
\begin{rem}\label{lem4}
Let $Q\in \mathcal{M}_{S,\boldsymbol{\Lambda }(\widetilde{\theta })}$. If
$P_{X_{1}}\neq Q_{X_{1}}$ or $P_{W_{1}}\neq Q_{W_{1}}$, applying
Lemma~\ref{lem3} together with the strong law of large numbers, it can be
easily seen that the probability measures\index{probability measure} $P$ and $Q$ are singular on
$\mathcal{F}_{\infty }$, implying that $P$ and $Q$ are equivalent of
$\mathcal{F}_{\infty }$ if and only if $P\upharpoonright \mathcal{F}_{\infty
}=Q\upharpoonright \mathcal{F}_{\infty }$ if and only if $P_{X_{1}}=
Q_{X_{1}}$ and $P_{W_{1}}= Q_{W_{1}}$.
\end{rem}

Before we formulate the inverse of Proposition~\ref{thm1} (i.e. that for
a given function $\gamma \in \mathcal{F}_{P,h}$ there exists a unique probability
measure\index{probability measure} $Q\in \mathcal{M}_{S,\boldsymbol{\Lambda }(\widetilde{\theta })}$ satisfying
\eqref{mart}) we remind a simple construction of \emph{canonical probability
spaces} admitting compound renewal processes.\index{compound renewal ! process}

By $(\varOmega \times \varXi ,\varSigma \otimes {H},P\otimes {R})$ we denote
the product probability space of the probability spaces $(\varOmega
,\varSigma ,P)$ and $(\varXi ,H,R)$. If $I$ is an arbitrary nonempty index
set, we write $P_{I}$ for the product measure on $\varOmega ^{I}$ and
$\varSigma _{I}$ for its domain.

\emph{{Throughout what follows, we put $\widetilde{\varOmega }:=\varUpsilon
^{\mathbb{N}}$, $\widetilde{\varSigma }:=\mathfrak{B}(\widetilde{\varOmega
})= \mathfrak{B}(\varUpsilon )_{\mathbb{N}}$, $\varOmega
:=\widetilde{\varOmega }\times \widetilde{\varOmega }$ and
$\varSigma :=\widetilde{\varSigma }\otimes \widetilde{\varSigma }$ for simplicity.}}%

For all $n\in \mathbb{N}$ and for any fixed $\theta \in D\subseteq
\mathbb{R}^{d}$, let $Q_{n}(\theta ):=\mathbf{{K}}(\theta )$ and $R_{n}:=R$
be probability measures\index{probability measure} on $\mathfrak{B}(\varUpsilon )$. Define the
probability measure\index{probability measure} $P$ on $\varSigma $ by means of $P:=\mathbf{K}(\theta
)_{\mathbb{N}}\otimes R_{\mathbb{N}}$, and for any $\omega =(w_{1},\ldots
,w_{n},\ldots ; x_{1},\ldots , x_{n},\ldots ) \in \varOmega $ put
$W_{n}(\omega ):=w_{n}$ and $X_{n}(\omega ):=x_{n}$. It then follows that
$X:=\{X_{n}\}_{n\in \mathbb{N}}$ is a claim size process satisfying the condition
$P_{X_{n}}=R$ for any $n\in \mathbb{N}$, and that $W:=\{W_{n}\}_{n\in
\mathbb{N}}$ is a $P$-independent claim interarrival process with
$P_{W_{n}}=\mathbf{K}(\theta )$ for any $n\in \mathbb{N}$. Putting
$T_{n}:=\sum_{k=1}^{n} W_{k}$ for any $n\in \mathbb{N}_{0}$ and
$T:=\{T_{n}\}_{n\in \mathbb{N}_{0}}$, we define by means of $N_{t}:=\sum_{n=1}^{\infty }\chi _{\{T_{n}\leq t\}}$ for any $t\geq 0$ the counting
process $N:=\{N_{t}\}_{t\in \mathbb{R}_{+}}$ induced by $T$ (cf. e.g.
\xch{\cite[Theorem 2.1.1]{Sc}}{\cite{Sc}, Theorem 2.1.1}). Setting $S_{t}:=\sum^{N_{t}}_{n=1} X_{n}$ for any
$t\geq 0$ and $S:=\{S_{t}\}_{t\in \mathbb{R}_{+}}$ we get that $S$ is a
$P$-CRP with parameters $\mathbf{K}(\theta )$ and $P_{X_{1}}$. Moreover,
according to Lemma~\ref{lem3} we have that $\varSigma =\mathcal{F}_{\infty
}^{(W,X)}=\mathcal{F}_{\infty }$.

The following proposition shows that after changing the measure the process
$S$ remains a compound renewal\index{compound renewal} one if the Radon--Nikod\'{y}m derivative has
the ``right'' structure on each $\sigma $-algebra $\mathcal{F}_{t}$. To
formulate it, we use the following notation and assumption.

%s3.1 #&#
\begin{symb}\label{symb3}
Let $\mathbf{K}(\theta )$ and $\boldsymbol{\Lambda }(\widetilde{\theta })$ be
probability distributions on $\mathfrak{B}(\varUpsilon )$ such that
$\mathbf{K}(\theta )\sim {\boldsymbol{\Lambda }}(\widetilde{\theta })$. For
any $n\in \mathbb{N}_{0}$ the class of all likelihood ratios
$g_{n}:=g_{\theta ,\widetilde{\theta },n}:\varUpsilon ^{n+1} \rightarrow
\varUpsilon $ defined by means of
\begin{equation*}
g_{n}(w_{1},\ldots ,w_{n},t):= \Biggl[\prod
_{j=1}^{n} \frac{d\boldsymbol{{\Lambda }}(\widetilde{\theta })}{d\mathbf{{K}}(\theta
)}(w_{j})
\Biggr]\cdot \frac{1-\boldsymbol{{\Lambda }}(\widetilde{\theta
})(t-w)}{1-\mathbf{{K}}(\theta )(t-w)}
\end{equation*}
for any $(w_{1},\ldots ,w_{n},t)\in \varUpsilon ^{n+1}$, where $w:=\sum^{n}_{j=1}w_{j}$, will be denoted by $\mathcal{G}_{n,\theta
,\widetilde{\theta }}$. Notation $\mathcal{G}_{\theta ,\widetilde{\theta }}$
stands for the set $\{g=\{g_{n}\}_{n\in \mathbb{N}_{0}} : g_{n}\in
\mathcal{G}_{n,\theta , \widetilde{\theta }}\text{ for any } n\in
\mathbb{N}_{0}\}$ of all sequences of elements of $\mathcal{G}_{n,\theta
,\widetilde{\theta }}$.
\end{symb}

\emph{Throughout what follows $\mathbf{K}(\theta )$, $\boldsymbol{\Lambda
}(\widetilde{\theta })$ and $g\in \mathcal{G}_{\theta ,\widetilde{\theta }}$
are as in Notation~\ref{symb3}, and $P$, $S$ are those constructed before
Notation~\ref{symb3}.}

%p3.1 #&#
\begin{prop}\label{thm2}
Let $\gamma \in \mathcal{F}_{P,h}$. Then for all $A\in \mathcal{F}_{t}$
the condition
\begin{equation*}
Q(A)=\int_{A} \Biggl[\prod_{j=1}^{N_{t}}
\,\bigl(h^{-1}\circ \gamma \circ X_{j}\bigr) \Biggr]\cdot
g_{N_{t}}(W_{1},\ldots ,W_{N_{t}},t)\,dP
\end{equation*}
determines a unique probability measure\index{probability measure} $Q\in
\mathcal{M}_{S,\boldsymbol{\Lambda }(\widetilde{\theta })}$.
\end{prop}

\begin{proof}
Fix an arbitrary $t\geq 0$, and define the set-functions
$\widecheck{Q}_{n}(\theta ), \widecheck{R}:\mathfrak{B}(\varUpsilon )
\rightarrow \mathbb{R}$ by means of $\widecheck{Q}_{n}(\theta
)(B_{1}):=\mathbb{E}_{P}[\chi _{W^{-1}_{1}(B_{1})} \cdot (
\frac{d\boldsymbol{\Lambda }(\widetilde{\theta })}{d\mathbf{K} (\theta )}
\circ W_{1})]$ and $\widecheck{R}(B_{2}):=\mathbb{E}_{P}[\chi
_{X^{-1}_{1}(B_{2})} \cdot (h^{-1} \circ \gamma \circ X_{1})]$ for any
$B_{1},B_{2}\in \mathfrak{B}(\varUpsilon )$, respectively. Applying a
monotone class\index{monotone class} argument it can be seen that $\widecheck{Q}_{n}(\theta
)=\boldsymbol{\Lambda } (\widetilde{\theta })$, while Lemma~\ref{lem2} (a)
$(i)$ implies that $\widecheck{R}$ is a probability measure.\index{probability measure} Therefore, we
may construct a probability measure\index{probability measure} $\widecheck{Q}:=\boldsymbol{\Lambda
}(\widetilde{\theta })_{\mathbb{N}} \otimes \widecheck{R}_{\mathbb{N}}$ on
$\varSigma $ such that $S$ is a $\widecheck{Q}$-CRP with parameters
$\boldsymbol{\Lambda }(\widetilde{\theta })$ and
$\widecheck{Q}_{X_{1}}=\widecheck{R}$, implying that
$\widecheck{Q}_{X_{1}}\sim P_{X_{1}}$ and $\widecheck{Q}_{W_{1}}\sim
P_{W_{1}}$. Applying now Proposition~\ref{thm1} we obtain
$\widecheck{Q}\upharpoonright \mathcal{F}_{t} \sim P\upharpoonright
\mathcal{F}_{t}$, or equivalently
\begin{equation*}
\widecheck{Q}(A)=\int_{A} \Biggl[\prod
_{j=1}^{N_{t}}\,\bigl(h^{-1}\circ \gamma
\circ X_{j}\bigr) \Biggr]\cdot g_{N_{t}}(W_{1},
\ldots ,W_{N_{t}},t)\,dP
\end{equation*}
for all $A\in \mathcal{F}_{t}$. Thus $Q\upharpoonright
\mathcal{F}_{t}=\widecheck{Q}\upharpoonright \mathcal{F}_{t}$; hence
$Q\upharpoonright \widecheck{\varSigma }=\widecheck{Q} \upharpoonright
\widecheck{\varSigma }$ where $\widecheck{\varSigma }:=\bigcup_{t\in
\mathbb{R}_{+}} \mathcal{F}_{t}$, implying that $Q$ is $\sigma $-additive on
$\widecheck{\varSigma }$ and that $\widecheck{Q}$ is the unique extension of
$Q$ on $\varSigma =\sigma (\widecheck{\varSigma })$.
\end{proof}

%r3.2 #&#
\begin{rem}\label{gen}
A well-known change of measure technique for compound renewal processes\index{compound renewal ! process} is to
markovize the process and then to change the measure (cf. e.g.
\xch{\cite[Chapter VI, Proposition 3.4]{asal}}{\cite{asal}, Chapter VI,
Proposition 3.4} or \xch{\cite[p.~139]{scm}}{\cite{scm}, p.~139}). Our method seems to simplify the
above one.
\end{rem}

The next result is the desired characterization. Its proof is an immediate
consequence of Propositions \ref{thm1} and \ref{thm2}.

%t3.1 #&#
\begin{thm}\label{thm!}
The following hold true:
\begin{enumerate}
\item[$(i)$] for any $Q\in \mathcal{M}_{S,\boldsymbol{\Lambda
    }(\widetilde{\theta })}$ there exists a $P_{X_{1}}$-a.s. unique
    function $\gamma \in \mathcal{F}_{P,h}$ satisfying condition
    \eqref{mart};
\item[$(ii)$] conversely, for any function $\gamma \in \mathcal{F}_{P,h}$
    there exists a unique probability measure\index{probability measure} $Q\in
    \mathcal{M}_{S,\boldsymbol{\Lambda }(\widetilde{\theta })}$ satisfying
    condition \eqref{mart}.
\end{enumerate}
\end{thm}

In order to formulate the next results of this section, let us denote by
$\widetilde{\mathcal{F}}_{P,\theta }$ the class of all real-valued
$\mathfrak{B}(\varUpsilon )$-measurable functions $\beta _{\theta }$, such
that $\beta _{\theta }:=\gamma +\alpha _{\theta }$, where $\gamma \in
\mathcal{F}_{P,\ln }$ and $\alpha _{\theta }$ is a real number depending
on~$\theta $.

The following result allows us to convert any compound renewal process\index{compound renewal ! process} into a
compound Poisson\index{compound Poisson} one through a change of measure.

%c3.1 #&#
\begin{cor}\label{cor4b}
If $W_{1}\in \mathcal{L}^{1}(P)$ then the following hold true:
\begin{enumerate}
\item[$(i)$] for any $\widetilde{\theta } \in \varUpsilon $ and any
    probability measure $Q\in \mathcal{M}_{S,\mathbf{Exp}(\widetilde{\theta
    })}$\index{probability measure} there exists a $P_{X_{1}}$-a.s. unique function $\beta _{\theta
    }\in \widetilde{\mathcal{F}}_{P,\theta }$ satisfying together with $Q$
    the conditions
%
%e{$\ast $} #&#
\begin{equation}
\label{ast} \alpha _{\theta }=\ln \widetilde{\theta }+\ln
\mathbb{E}_{P}[W_{1}] \tag{\texorxml{$\ast$}{\textasteriskcentered}}
\end{equation}
and
%
%e{RPM} #&#
\begin{equation}
\label{martPPb}\tag{RPM} Q(A)=\int_{A}
M^{(\beta )}_{t}(\theta )\,dP\quad \text{for all } A \in
\mathcal{F}_{t},
\end{equation}
where $M^{(\beta )}_{t}(\theta ):= \frac{e^{\sum _{j=1}^{N_{t}} \beta
_{\theta }(X_{j})-\widetilde{\theta }\cdot (t-T_{N_{t}})}\cdot
(\widetilde{\theta }\cdot \mathbb{E}_{P}[W_{1}]
)^{-N_{t}}}{1-\mathbf{{K}}(\theta )(t-T_{N_{t}})} \cdot   [\prod_{j=1}^{N_{t}}\,\frac{dQ_{W_{1}}}{dP_{W_{1}}}(W_{j})   ]$;
\item[$(ii)$] conversely, for any function $\beta _{\theta }\in
    \widetilde{\mathcal{F}}_{P,\theta }$ there exist a $\widetilde{\theta
    }\in \varUpsilon $ and a unique probability measure\index{probability measure} $Q\in
    \mathcal{M}_{S,\mathbf{Exp}(\widetilde{\theta })}$ satisfying together
    with $\beta _{\theta }$ conditions \eqref{ast} and \eqref{martPPb}.
\end{enumerate}
\end{cor}

\begin{proof}
Fix an arbitrary $t\geq 0$.

For $(i)$: Under the assumptions of statement $(i)$, according to
Theorem~\ref{thm!} $(i)$ there exists a $P_{X_{1}}$-a.s. unique function
$\gamma \in \mathcal{F}_{P,\ln }$ defined by means of $\gamma :=\ln f$, where
$f$ is a Radon--Nikod\'{y}m derivative of $Q_{X_{1}}$ with respect to
$P_{X_{1}}$, such that
%
%e5 #&#
\begin{equation}
Q(A)=\int_{A} \frac{e^{\sum _{j=1}^{N_{t}} \gamma (X_{j}) -\widetilde{\theta
}\cdot (t-T_{N_{t}})}}{ 1-\mathbf{{K}}(\theta )(t-T_{N_{t}})} \cdot \Biggl[\prod
_{j=1}^{N_{t}}\,\frac{dQ_{W_{1}}}{dP_{W_{1}}}(W_{j})
\Biggr]\,dP \label{1111c}
\end{equation}
for all $A\in \mathcal{F}_{t}$. Define $\alpha _{\theta }:=\ln
\widetilde{\theta }+\ln \mathbb{E}_{P}[W_{1}]$, and put $\beta _{\theta
}:=\gamma +\alpha _{\theta }$. It then follows that $\beta _{\theta }\in
\widetilde{\mathcal{F}}_{P,\theta }$ and that condition \eqref{ast} is valid.
The latter together with condition \eqref{1111c} implies condition
\eqref{martPPb}.

\xch{For}{ad} $(ii)$: Let $\beta _{\theta }=\gamma +\alpha _{\theta }\in
\widetilde{\mathcal{F}}_{P, \theta }$ and define $\widetilde{\theta }:=
\frac{e^{\alpha _{\theta }}}{\mathbb{E}_{P}[W_{1}]}$. By Theorem~\ref{thm!}
$(ii)$ for the function $\gamma =\beta _{\theta }-\alpha _{\theta }$ there
exists a unique probability measure\index{probability measure} $Q\in
\mathcal{M}_{S,\mathbf{Exp}(\widetilde{\theta })}$ satisfying condition
\eqref{mart} or equivalently condition \eqref{martPPb}.
\end{proof}

%r3.3 #&#
\begin{rem}\label{cla}
{\textbf{(a)}} In the special case $P\in \mathcal{M}_{S,\mathbf{Exp}(\theta
)}$, Corollary~\ref{cor4b} yields the main result of Delbaen and Haezendonck
\xch{\cite[Proposition 2.2]{dh}}{\cite{dh}, Proposition 2.2}.

\textbf{(b)} Theorem~\ref{thm!} remains true if we replace the classes
$\mathcal{M}_{S,\boldsymbol{\Lambda }(\widetilde{\theta })}$ and
$\mathcal{F}_{P,h}$ by their subclasses $\mathcal{M}^{\ell
}_{S,\boldsymbol{\Lambda }(\widetilde{\theta })}:= \{Q\in
\mathcal{M}_{S,\boldsymbol{\Lambda }(\widetilde{\theta })} :
\mathbb{E}_{Q}[X^{\ell }_{1}]<\infty ]\}$ and $\mathcal{F}^{\ell
}_{P,h}:=\{\gamma \in \mathcal{F}_{P,h} : \mathbb{E}_{P}[X^{\ell }_{1}\cdot
(h^{-1}\circ \gamma \circ X_{1})]<\infty \}$ for $\ell =1,2$, respectively.
As a consequence, Corollary~\ref{cor4b} remains true if we replace the class
$\widetilde{\mathcal{F}}_{P,\theta }$ by its subclass
$\widetilde{\mathcal{F}}^{\ell }_{P,\theta }:=\{\beta _{\theta }= \gamma +
\alpha _{\theta }: \gamma \in \mathcal{F}^{\ell }_{P,\ln }\text{ and }\alpha
_{\theta }\in \mathbb{R}\}$ for $\ell =1,2$.
\end{rem}

The following example translates the results of Corollary~\ref{cor4b} to a
well-known compound renewal process\index{compound renewal ! process} appearing in applications.

%e3.1 #&#
\begin{ex}\label{ex7}
Fix an arbitrary $t\geq 0$, let $\theta :=(\xi ,2)\in D:=\varUpsilon ^{2}$,
and let $P\in \mathcal{M}_{S,\mathbf{Ga}(\theta )}$ such that
$P_{X_{1}}=\mathbf{Ga}(\eta )$, where $\eta :=(b,2)\in D$. Let
$\widetilde{\theta }\in \varUpsilon $ and $Q\in
\mathcal{M}_{S,\mathbf{Exp}(\widetilde{\theta })}$ such that
$Q_{X_{1}}=\mathbf{Exp}(\zeta )$, where $\zeta $ is a positive real constant.
By Corollary~\ref{cor4b} $(i)$, there exists a $P_{X_{1}}$-a.s. unique
function $\beta _{\theta }:=\gamma +\alpha _{\theta }\in
\widetilde{\mathcal{F}}_{P, \theta }$, where $\gamma (x):=\ln \frac{\zeta
\cdot e^{-\zeta \cdot x}}{b^{2}\cdot x\cdot e^{-b\cdot x}}$ for any $x\in
\varUpsilon $ and $\alpha _{\theta }:=\ln \widetilde{\theta }+\ln
\mathbb{E}_{P}[W_{1}]= \ln \frac{2\cdot \widetilde{\theta }}{\xi }$,
satisfying together with $Q$ the condition
%
%e6 #&#
\begin{equation}
Q(A)=\int_{A}\, \biggl(\frac{1}{2\xi }
\biggr)^{N_{t}}\cdot \frac{e^{\sum
_{j=1}^{N_{t}} \beta _{\theta }(X_{j})-t\cdot \widetilde{\theta }+t \xi
}}{[\prod^{N_{t}}_{j=1}W_{j}]\cdot   (1+ \xi \cdot (t-T_{N_{t}})  )} \,dP \label{ex12}
\end{equation}
for all $A\in \mathcal{F}_{t}$.

Conversely, let $\zeta $ be as above and consider the function $\beta
_{\theta }:=\gamma +\alpha _{\theta }$,\break where $\gamma (x):=\ln \frac{\zeta
\cdot e^{-\zeta \cdot x}}{b^{2}\cdot x\cdot e^{-b\cdot x}}$ for any $x\in
\varUpsilon $ and $\alpha _{\theta }\in \mathbb{R}$. It then follows easily
that $\mathbb{E}_{P}[e^{\gamma (X_{1})}]=1$, implying that $\gamma \in
\mathcal{F}_{P,\ln }$; hence $\beta _{\theta }\in
\widetilde{\mathcal{F}}_{P,\theta }$. Thus, we may apply
Corollary~\ref{cor4b} $(ii)$ to get a $\widetilde{\theta }\in \varUpsilon $
and a unique probability measure\index{probability measure} $Q\in
\mathcal{M}_{S,\mathbf{Exp}(\widetilde{\theta })}$ satisfying together with
$\beta _{\theta }$ conditions \eqref{ast} and \eqref{ex12}. But then applying
Lemma~\ref{lem2} (a), we get
\begin{equation*}
Q_{X_{1}}(B)=\mathbb{E}_{P}\bigl[\chi _{X^{-1}_{1}(B)}
\cdot e^{\gamma (X_{1})}\bigr]= \int_{B} \zeta \cdot
e^{-\zeta \cdot x}\,\lambda (dx)\quad\text{for any } B\in \mathfrak{B}(
\varUpsilon ),
\end{equation*}
implying that $Q_{X_{1}}=\mathbf{Exp}(\zeta )$.
\end{ex}

%s4 #&#
\section{Applications}\label{CRPM}

In this section we first show that a martingale approach to premium
calculation principles leads\index{premium calculation principles} in the case of CRPs to CPPs, providing in this
way a method to find progressively equivalent martingale measures.\index{equivalent martingale measures} Next,
using our results we show that if $\widetilde{\mathcal{F}}^{2}_{P,\theta
}\neq \emptyset $ then there exist \emph{canonical} price processes (called
\emph{claim surplus processes} in Risk Theory) satisfying the condition of no
free lunch with vanishing risk.

In order to present the results of this section we recall the following
notions. For a given real-valued process $Y:=\{Y_{t}\}_{t\in \mathbb{R}_{+}}$
on $(\varOmega ,\varSigma )$ a probability measure\index{probability measure} $Q$ on $\varSigma $ is
called a \textbf{martingale measure\index{martingale measure}} for $Y$, if $Y$ is a $Q$-martingale. We
will say that $Y$ satisfies condition (PEMM) if there exists a
\textbf{progressively equivalent martingale measure\index{equivalent martingale measures}} (PEMM for short) for
$Y$, i.e. a probability measure\index{probability measure} $Q$ on $\varSigma $ such that
$Q\upharpoonright \mathcal{F}_{t}\sim P\upharpoonright \mathcal{F}_{t}$ for
any $t\geq 0$ and $Y$ is a $Q$-martingale. Moreover, let $T>0$,
$\mathbb{T}:=[0,T]$, $Q_{T}:=Q\upharpoonright \mathcal{F}_{T}$,
$Y_{\mathbb{T}}:=\{Y_{t}\}_{t\in \mathbb{T}}$ and
$\mathcal{F}_{\mathbb{T}}:=\{\mathcal{F}_{t}\}_{t\in \mathbb{T}}$. We will
say that the process $Y_{\mathbb{T}}$ satisfies condition (EMM) if there
exists an \textbf{equivalent martingale measures\index{equivalent martingale measures}} for $Y_{\mathbb{T}}$, i.e.
a probability measure\index{probability measure} $Q_{T}$ on $\mathcal{F}_{T}$ such that $Q_{T}\sim
P_{T}$ and $Y_{\mathbb{T}}$ is a
$(Q_{T},\mathcal{F}_{\mathbb{T}})$-martingale.

Suppose that $X_{1},W_{1}\in \mathcal{L}^{1}(P)$ and define the
\textbf{premium density\index{premium density}} as
\begin{equation*}
p(P):=\frac{\mathbb{E}_{P}[X_{1}]}{\mathbb{E}_{P}[W_{1}]}\in \varUpsilon .
\end{equation*}
Consider the process $Z(P):=\{Z_{t}\}_{t\in \mathbb{R}_{+}}$ with
$Z_{t}:=S_{t}-t\cdot p(P)$ for any $t\geq 0$. The following auxiliary result
could be of independent interest, since it says that if $S$ is under $P$ a
CRP and the process $Z(P)$ is a $P$-martingale, then $N_{t}$ must have a
Poisson distribution so that $S$ is actually a CPP.

%p4.1 #&#
\begin{prop}\label{thm4}
Let $\widetilde{\theta }:=\frac{1}{\mathbb{E}_{P}[W_{1}]}$. Consider the
following statements:
\begin{enumerate}
\item[$(i)$] $P$ is a martingale measure\index{martingale measure} for $Z(P)$;
\item[$(ii)$] $P\in \mathcal{M}^{1}_{S,\mathbf{Exp}(\widetilde{\theta })}$;
\item[$(iii)$] $P$ is a martingale measure\index{martingale measure} for $Z(P)$ such that $Z_{t}\in
    \mathcal{L}^{2}(P)$ for any $t\geq 0$;
\item[$(iv)$] $P\in \mathcal{M}^{2}_{S,\mathbf{Exp}(\widetilde{\theta })}$.
\end{enumerate}
Then statements $(i)$ and $(ii)$ as well as statements $(iii)$ and $(iv)$ are
equivalent. Moreover, if $X_{1}\in \mathcal{L}^{2}(P)$ then all statements
$(i)$ to $(iv)$ are equivalent.
\end{prop}

\begin{proof}
Fix an arbitrary $t\geq 0$.

For $(i)\Longrightarrow (ii)$: Since $P$ is martingale measure\index{martingale measure} for $Z(P)$, we
have $\mathbb{E}_{P}[Z_{t}]=0$ implying $\mathbb{E}_{P}[S_{t}]=t\cdot
\frac{\mathbb{E}_{P}[X_{1}]}{\mathbb{E}_{P}[W_{1}]}$, or equivalently
$\mathbb{E}_{P}[N_{t}]=t\cdot \widetilde{\theta }$.

%i4 #&#

\begin{Claim}
The following are equivalent:
\begin{itemize}
\item[(a)] $N$ is a $P$-Poisson process with parameter $\theta $;
\item[(b)] $\mathbb{E}_{P}[N_{t}] = t\theta $.
\end{itemize}
\end{Claim}

\begin{proof}
The above claim is well known (cf. e.g. \xch{\cite[Remark 21, p.~110]{se}}{\cite{se}, Remark 21 p.~110}), but
since we have not seen its proof anywhere, we insert it for completeness.
The implication $(a)\Longrightarrow (b)$ is immediate.

For $(b)\Longrightarrow (a)$: To prove this implication, let us recall that
the \textbf{renewal function associated with the distribution}
$\mathbf{K}(\theta )$ is defined by
\begin{equation*}
U(u):=\sum^{\infty }_{n=0} {
\mathbf{K}}^{\ast n}(\theta ) (u)\quad \text{for any }u\in \mathbb{R}
\end{equation*}
where $\mathbf{K}^{\ast n}(\theta )$ is the $n$-fold convolution of
$\mathbf{K}(\theta )$ (cf. e.g. \xch{\cite[Definition 17, p.~108]{se}}{\cite{se}, Definition 17, p.~108}). Clearly
$U(u)=1+ \mathbb{E}_{P}[N_{u}]$ for any $u\geq 0$. Assuming that
$\mathbb{E}_{P}[N_{t}]=t\theta $, we get $U(t)=1+t\theta $, implying that the
Laplace--Stieltjes transform $\widehat{U}(s)$ of $U(t)$ is given by
\begin{equation*}
\widehat{U}(s)=\int_{\mathbb{R}_{+}} e^{-s\cdot u}
\,dU(u)=e^{-s \cdot
0}\cdot U(0)+\int_{0}^{\infty }
\theta e^{-s\cdot u}\,du= \frac{s+\theta
}{s}\quad \text{for every } s\geq
0,
\end{equation*}
where the second equality follows from the fact that $\int_{\mathbb{R}_{+}}
e^{-s\cdot u}\,dU(u)$ is a Riemann--Stieltjes integral and $U$ has a density
for $u>0$, $U(u)=0$ for $u<0$ and it has a unit jump at $u=0$ (cf. e.g.
\xch{\cite[pp.~108--109]{se}}{\cite{se}, pp.~108--109}). It then follows that
\begin{equation*}
\widehat{\mathbf{K}}(\theta ) (s)= \frac{\widehat{U}(s)-1}{\widehat{U}(s)}=\frac{\theta }{\theta +s}
\quad \text{for any } s\geq 0,
\end{equation*}
where $\widehat{\mathbf{K}}(\theta )$ denotes the Laplace--Stieltjes transform
of the distribution of $W_{n}$ for any $n\in \mathbb{N}$ (cf. \xch{e.g.}{e.g}
\xch{\cite[Proposition 20, p.~109]{se}}{\cite{se}, Proposition 20, p.~109}); hence $P_{W_{n}}=\mathbf{Exp}(\theta )$
for any $n\in \mathbb{N}$. But since $W$ is also $P$-independent, it follows
that $N$ is a $P$-Poisson process with parameter $\theta $ (cf. e.g.
\xch{\cite[Theorem 2.3.4]{Sc}}{\cite{Sc}, Theorem 2.3.4}).
\end{proof}

Thus, according to the above claim statement $(ii)$ follows.

For $(ii)\Longrightarrow (i)$: Since $P\in
\mathcal{M}^{1}_{S,\mathbf{Exp}(\widetilde{\theta })}$, it follows that $S$
has independent increments (cf. e.g. \xch{\cite[Theorem 5.1.3]{Sc}}{\cite{Sc}, Theorem 5.1.3}). Thus, for all
$u\in [0,t]$ and $A\in \mathcal{F}_{u}$ we get
\begin{align*}
&\int_{A} \bigl(S_{t}-\mathbb{E}_{P}[S_{t}]
\bigr)-\bigl(S_{u}-\mathbb{E}_{P}[S_{u}]
\bigr)\,dP\\
&\quad= \int_{\varOmega }\chi _{A}\,dP\cdot \int
_{\varOmega }\bigl((S_{t}-S_{u})-
\mathbb{E}_{P}[S_{t}-S_{u}]\bigr)\,dP=0,
\end{align*}
implying that the process $\{S_{t}-\mathbb{E}_{P}[S_{t}]\}_{t\in
\mathbb{R}_{+}}$ is a $P$-martingale. But since $\mathbb{E}_{P}[S_{t}]=t\cdot
\mathbb{E}_{P}[S_{1}]$, statement $(i)$ follows.

For $(iii)\Longrightarrow (iv)$: Since $P$ is a martingale measure\index{martingale measure} for $Z(P)$,
it follows by the equivalence of statements $(i)$ and $(ii)$ that $P\in
\mathcal{M}^{1}_{S,\mathbf{Exp}(\widetilde{\theta })}$. But since $Z_{t}\in
\mathcal{L}^{2}(P)$, we have $\operatorname{Var}_{P}[Z_{t}]=\mathbb{E}_{P}[N_{t}]\cdot
\operatorname{Var}_{P}[X_{1}]+\operatorname{Var}_{P}[N_{t}] \cdot \mathbb{E}^{2}_{P}[X_{1}]<\infty $, where
$\operatorname{Var}_{P}$ denotes the variance under the measure $P$; hence
$\operatorname{Var}_{P}[X_{1}]<\infty $, implying statement $(iv)$.

For $(iv)\Longrightarrow (iii)$: Since $P\in
\mathcal{M}^{2}_{S,\mathbf{Exp}(\widetilde{\theta })}$ and
$\mathcal{M}^{2}_{S,\mathbf{Exp}(\widetilde{\theta })}\subseteq
\mathcal{M}^{1}_{S,\mathbf{Exp}(\widetilde{\theta })}$, it follows again by
the equivalence of statements $(i)$ and $(ii)$ that $P$ is a martingale
measure\index{martingale measure} for $V$. But
$\operatorname{Var}_{P}[Z_{t}]=\operatorname{Var}_{P}[S_{t}]=\mathbb{E}_{P}[N_{t}]\cdot
\operatorname{Var}_{P}[X_{1}]+\operatorname{Var}_{P}[N_{t}] \cdot \mathbb{E}^{2}_{P}[X_{1}]<\infty $, where
the inequality follows by the fact that $P\in
\mathcal{M}^{2}_{S,\mathbf{Exp}(\widetilde{\theta })}$; hence statement
$(iii)$ follows.

Moreover, assuming statement $(ii)$ and $X_{1}\in \mathcal{L}^{2}(P)$, we get
immediately state\-ment $(iv)$, implying that all statements $(i)$--$(iv)$ are
equivalent.
\end{proof}

In the next proposition we find out a wide class of \emph{canonical}
processes convert\-ing the progressively equivalent measures $Q$ of
Theorem~\ref{thm!} into martingale mea\-sures.\index{martingale measure} In this way, a
characterization of all progressively equivalent martingale measures,\index{equivalent martingale measures} similar
to that of Theorem~\ref{thm!}, is provided.

%p4.2 #&#
\begin{prop}\label{mar}
If $\ell =1,2$ and $P\in \mathcal{M}^{\ell }_{S,\mathbf{K}(\theta )}$ the
following hold true:
\begin{enumerate}
\item[$(i)$] for every $\widetilde{\theta }\in \varUpsilon $ and $Q\in
    \mathcal{M}^{\ell }_{S,\mathbf{Exp}(\widetilde{\theta })}$ there exists
    a $P_{X_{1}}$-a.s. unique function $\beta _{\theta }\in
    \widetilde{\mathcal{F}}^{\ell }_{P,\theta }$ satisfying together with
    $Q$ conditions \eqref{ast} and \eqref{martPPb}, and the process
    $V:=\{V_{t}\}_{t\in \mathbb{R}_{+}}$, defined by means of
    $V_{t}:=S_{t}-t\cdot \frac{\mathbb{E}_{P}[X_{1}\cdot e^{\beta _{\theta
    }(X_{1})}]}{\mathbb{E}_{P}[W_{1}]}$ for any $t\geq 0$, such that $Q$ is
    a PEMM for $V$;
\item[$(ii)$] conversely, for every function $\beta _{\theta }\in
    \widetilde{\mathcal{F}}^{\ell }_{P,\theta }$ and for the process $V$
    defined in $(i)$, there exist a $\widetilde{\theta }\in \varUpsilon $
    and a unique probability measure\index{probability measure} $Q\in \mathcal{M}^{\ell
    }_{S,\mathbf{Exp}(\widetilde{\theta })}$ satisfying together with
    $\beta _{\theta }$ conditions \eqref{ast} and \eqref{martPPb}, and such
    that $Q$ is a PEMM for $V$.
\end{enumerate}
In both cases $V=Z(Q)$.
\end{prop}

\begin{proof}
Fix $\ell =1$ or $\ell =2$.

For $(i)$: Under the assumptions of $(i)$, by Corollary~\ref{cor4b} $(i)$ and
Remark~\ref{cla} (b) there exists a $P_{X_{1}}$-a.s. unique function $\beta
_{\theta }\in \widetilde{\mathcal{F}}^{\ell }_{P,\theta }$ satisfying
together with $Q$ conditions \eqref{ast} and \eqref{martPPb}. It then follows by
Lemma~\ref{lem2} (a) and condition \eqref{ast} that $V=Z(Q)$; hence by
Proposition~\ref{thm4} we get that $Q$ is a PEMM for $V$.

For $(ii)$: Under the assumptions of $(ii)$, by Corollary~\ref{cor4b} $(ii)$
and Remark~\ref{cla} (b) there exist a $\widetilde{\theta }\in \varUpsilon $
and a unique probability measure\index{probability measure} $Q\in \mathcal{M}^{\ell
}_{S,\mathbf{Exp}(\widetilde{\theta })}$ satisfying together with $\beta
_{\theta }$ conditions \eqref{ast} and \eqref{martPPb}; hence according to
Proposition~\ref{thm4} the process $Z(Q)$ is a $Q$-martingale. Again by
Lemma~\ref{lem2} (a) and condition \eqref{ast} we obtain that $V=Z(Q)$.
\end{proof}

The next theorem connects our results with the basic notion of no free lunch
with vanishing risk ((NFLVR) for short) (see \xch{\cite[Definition 8.1.2]{ds}}{\cite{ds}, Definition 8.1.2}) of
Mathematical Finance.

%t4.1 #&#
\begin{thm}\label{nfl}
Let $P\in \mathcal{M}^{2}_{S,\mathbf{K}(\theta )}$, $\beta _{\theta }\in
\widetilde{\mathcal{F}}^{2}_{P,\theta }$ and $V$ be as above. There exist a
$\widetilde{\theta }\in \varUpsilon $ and a unique probability measure\index{probability measure} $Q\in
\mathcal{M}^{2}_{S,\mathbf{Exp}(\widetilde{\theta })}$ satisfying together
with $\beta _{\theta }$ conditions \eqref{ast} and \eqref{martPPb}, and such
that for every $T>0$ the process $V_{\mathbb{T}}:=\{V_{t}\}_{t\in
\mathbb{T}}$ satisfies condition (NFLVR).
\end{thm}

\begin{proof}
Fix an arbitrary $T>0$ and let $\beta _{\theta }\in
\widetilde{\mathcal{F}}^{2}_{P,\theta }$. By Proposition~\ref{mar} $(ii)$
there exist a $\widetilde{\theta }\in \varUpsilon $ and a unique probability
measure\index{probability measure} $Q\in \mathcal{M}^{2}_{S,\mathbf{Exp}(\widetilde{\theta })}$
satisfying together with $\beta _{\theta }$ conditions \eqref{ast} and
\eqref{martPPb}, and such that $V$ is a $Q$-martingale with $V_{t}\in
\mathcal{L}^{2}(Q)$ for any $t\geq 0$; hence $V_{\mathbb{T}}$ is a
$(Q_{T},\mathcal{F}_{\mathbb{T}})$-martingale, implying that it is a
$(Q_{T},\mathcal{F}_{\mathbb{T}})$-semi-martingale (cf. e.g.
\xch{\cite[Definition 7.1.1]{vw}}{\cite{vw}, Definition 7.1.1}). The latter
implies that $V_{\mathbb{T}}$ is also a
$(P_{T},\mathcal{F}_{\mathbb{T}})$-semi-martingale since $Q_{T}\sim P_{T}$
(cf. e.g. \xch{\cite[Theorem 10.1.8]{vw}}{\cite{vw}, Theorem 10.1.8}). But since the process $V$ satisfies
condition (PEMM) we have that $V_{\mathbb{T}}$ satisfies condition (EMM).
Thus, applying the Fundamental Theorem of Asset Pricing (FTAP for short) for
unbounded stochastic processes, see \xch{\cite[Theorem 14.1.1]{ds}}{\cite{ds}, Theorem 14.1.1}, we obtain that
the process $V_{\mathbb{T}}$ satisfies condition (NFLVR).
\end{proof}

%r4.1 #&#
\begin{rem}\label{ush}
It is well known that the FTAP of Delbaen and Schachermayer uses P.A.~Meyer's
\textit{usual conditions} (cf. e.g. \xch{\cite[Definition 2.1.5]{vw}}{\cite{vw}, Definition 2.1.5}). These
conditions play a fundamental role in the definition of the stochastic
integral with respect to a (\mbox{semi-})mar\-tingale. Nevertheless, the stochastic
integral can be defined for any semi-martingale without the usual conditions
(see \xch{\cite[pp.~22--23 and p.~150]{vw}}{\cite{vw}, pp.~22--23 and p.~150}). As a consequence, the \textit{easy}
implication of the FTAP of Delbaen and Schachermayer (i.e. (EMM)
$\Longrightarrow $ (NFLVR)) holds true without the usual conditions.
\end{rem}

We have seen that the initial probability measure $P$\index{probability measure} can be replaced by
another progressively equivalent probability measure\index{probability measure} $Q$ such that $S$ is
converted into a $Q$-CPP. The idea is to define a probability measure\index{probability measure} $Q$ in
order to give more weight to less favourable events. More precisely $Q$ must
be defined in such a way that the corresponding premium density $p(Q)$
includes the \emph{safety loading}, i.e. $p(P)<p(Q)$. This led Delbaen and
Haezendonck to define a \textbf{premium calculation principle\index{premium calculation principles}} as a
probability measure\index{probability measure} $Q\in \mathcal{M}^{1}_{S,\mathbf{Exp}(\lambda )}$, for
some $\lambda \in \varUpsilon $ (compare \xch{\cite[Definition 3.1]{dh}}{\cite{dh}, Definition 3.1}).

In the next Examples \ref{expcp0} to \ref{expcp7}, applying
Proposition~\ref{mar} and Theorem~\ref{nfl}, we show how to construct premium
calculation principles\index{premium calculation principles} $Q$ satisfying the desired property
$p(P)<p(Q)<\infty$, and such that for any $T>0$ the process $V_{\mathbb{T}}$
has the property of (NFLVR). For a discussion on how to rediscover some
well-known premium calculation principles\index{premium calculation principles} in the frame of classical Risk Theory
using change of measures techniques we refer to \xch{\cite[Examples 3.1 to~3.3]{dh}}{\cite{dh}, Examples 3.1
to~3.3}.

%e4.1 #&#
\begin{ex}\label{expcp0}
Let $\theta :=(\xi ,k)\in D:=\varUpsilon ^{2}$, and let $P\in
\mathcal{M}^{2}_{S,\mathbf{Ga}(\theta )}$ be such that
$P_{X_{1}}=\mathbf{Ga}(\eta )$, where $\eta :=(\zeta ,2)\in D$. Consider the
real-valued function $\beta _{\theta }:=\gamma +\alpha _{\theta }$ with
$\gamma (x):=\ln \frac{\mathbb{E}_{P}[X_{1}]}{2c}-\ln x+
\frac{2(c-1)}{c\mathbb{E}_{P}[X_{1}]}\cdot x$ for any $x\in \varUpsilon $,
where $c>2$ is a real constant, and $\alpha _{\theta }:=\ln (\frac{\xi
}{d}\cdot \mathbb{E}_{P}[W_{1}])$, where $d<k$ is a positive constant. It can
be easily seen that $\mathbb{E}_{P}[e^{\gamma (X_{1})}]=1$ and
$\mathbb{E}_{P}[X_{1}^{2}\cdot e^{\gamma (X_{1})}]= \frac{2c^{2}}{\zeta
}<\infty $, implying $\gamma \in \mathcal{F}^{2}_{P,\ln }$; hence $\beta
_{\theta }\in \widetilde{\mathcal{F}}^{2}_{P,\theta }$. Define
$\widetilde{\theta }$ by means of $\widetilde{\theta }:= e^{\alpha _{\theta
}}/\mathbb{E}_{P}[W_{1}]$. Thus, due to Proposition~\ref{mar} $(ii)$, there
exists a unique premium calculation principle\index{premium calculation principles} $Q\in
\mathcal{M}^{2}_{S,\mathbf{Exp}(\widetilde{\theta })}$ satisfying conditions
\eqref{ast} and \eqref{martPPb}, and such that $Q$ is a PEMM for the process
$V$ with $V_{t}:=S_{t}-t\cdot \frac{\xi }{d}\cdot
\frac{\mathbb{E}_{P}[X_{1}]}{2c}\cdot \mathbb{E}_{P}  [e^{ \frac{2\cdot
(c-1)}{c\cdot \mathbb{E}_{P}[X_{1}]}\cdot X_{1}}  ] \in
\mathcal{L}^{2}(Q)$ for any $t\geq 0$. Therefore, applying Lemma~\ref{lem2}
(a) we get
\begin{equation*}
Q_{X_{1}}(B)=\int_{B} \frac{\zeta }{c}
\cdot e^{-\frac{\zeta }{c} \cdot
x}\,\lambda (dx)\quad\text{for any } B\in \mathfrak{B}(
\varUpsilon ),
\end{equation*}
implying that $Q_{X_{1}}=\mathbf{Exp}(\frac{\zeta }{c})$; hence
$p(P)<p(Q)<\infty $. In particular, according to Theorem~\ref{nfl} for any
$T>0$ the process $V_{\mathbb{T}}$ satisfies the (NFLVR) condition.
\end{ex}

%e4.2 #&#
\begin{ex}\label{expcp6b}
Let $\theta :=(k,b)\in D:=\varUpsilon ^{2}$, let $\mathbf{W}(\theta )$ be the
Weibull distribution over $\mathfrak{B}(\varUpsilon )$ defined by %means of
\begin{equation*}
\mathbf{W}(\theta ) (B):=\int_{B} \frac{k}{b^{k}}
\cdot x^{k-1}\cdot e^{-(x/b)^{k}} \,\lambda (dx)\quad
\text{for any } B\in \mathfrak{B}(\varUpsilon ),
\end{equation*}
and let $P\in \mathcal{M}^{2}_{S,\mathbf{W}(\theta )}$ such that
$P_{X_{1}}=\mathbf{Exp}(\eta )$, where $\eta \in \varUpsilon $. Consider the
real-valued function $\beta _{\theta }:=\gamma +\alpha _{\theta }$ with
$\gamma (x):=\ln (1-c\cdot \mathbb{E}_{P}[X_{1}])+c\cdot x$ for any $x\in
\varUpsilon $, where $c<\eta $ is a positive constant, and $\alpha _{\theta
}:=0$. It can be easily seen that $\mathbb{E}_{P}[e^{\gamma (X_{1})}]=1$ and
$\mathbb{E}_{P}[X_{1}^{2}\cdot e^{\gamma (X_{1})}]= \frac{2}{(\eta
-c)^{2}}<\infty $, implying $\gamma \in \mathcal{F}^{2}_{P,\ln }$; hence
$\beta _{\theta }\in \widetilde{\mathcal{F}}^{2}_{P,\theta }$. Define
$\widetilde{\theta }$ by %means of
$\widetilde{\theta }:= e^{\alpha _{\theta
}}/\mathbb{E}_{P}[W_{1}]$. Applying now Proposition~\ref{mar} $(ii)$ we get
that there exists a unique premium calculation principle\index{premium calculation principles} $Q\in
\mathcal{M}^{2}_{S,\mathbf{Exp}(\widetilde{\theta })}$ satisfying conditions
\eqref{ast} and \eqref{martPPb}, and such that $Q$ is a PEMM for the process
$V$ with $V_{t}:=S_{t}-t\cdot \frac{(1-c\cdot \mathbb{E}_{P}[X_{1}])\cdot
\mathbb{E}_{P}[X_{1}\cdot e^{c\cdot X_{1}}]}{b\cdot \Gamma (1+1/k)} \in
\mathcal{L}^{2}(Q)$ for any $t\geq 0$. The latter together with
Lemma~\ref{lem2} (a) yields
\begin{equation*}
Q_{X_{1}}(B)=\int_{B} (\eta -c)\cdot
e^{-(\eta -c)\cdot x}\,\lambda (dx) \quad\text{for any } B\in \mathfrak{B}(
\varUpsilon ),
\end{equation*}
implying that $Q_{X_{1}}=\mathbf{Exp}(\eta -c)$. Thus, $p(P)<p(Q)<\infty $.
In particular, according to Theorem~\ref{nfl} for any $T>0$ the process
$V_{\mathbb{T}}$ satisfies the (NFLVR) condition.
\end{ex}

In our next example we show how one can obtain the Esscher principle by
applying Proposition~\ref{mar} $(ii)$.

%e4.3 #&#
\begin{ex}\label{expcp7}
Take $\theta :=(\xi ,2)\in D:=\varUpsilon ^{2}$, and let $P\in
\mathcal{M}^{2}_{S,\mathbf{Ga}(\theta )}$ such that
$P_{X_{1}}=\mathbf{Ga}(\eta )$, where $\eta :=(b,a)\in D$. Consider the
real-valued function $\beta _{\theta }:=\gamma +\alpha _{\theta }$ with
$\gamma (x):=c\cdot x-\ln \mathbb{E}_{P}[e^{c\cdot X_{1}}]$ for any $x\in
\varUpsilon $, where $c<b$ is a positive constant, and $\alpha _{\theta
}:=0$. It can be easily seen that $\mathbb{E}_{P}[e^{\gamma (X_{1})}]=1$ and
$\mathbb{E}_{P}[X_{1}^{2}\cdot e^{\gamma (X_{1})}]= \frac{a\cdot
(a+1)}{(b-c)^{2}}<\nobreak\infty $, implying $\gamma \in \mathcal{F}^{2}_{P,\ln }$;
hence $\beta _{\theta }\in \widetilde{\mathcal{F}}^{2}_{P,\theta }$.
Define $\widetilde{\theta }$ by %means of
$\widetilde{\theta }:= e^{\alpha _{\theta
}}/\mathbb{E}_{P}[W_{1}]$. Thus, due to Proposition~\ref{mar} $(ii)$ there
exists a unique premium calculation principle\index{premium calculation principles} $Q\in
\mathcal{M}^{2}_{S,\mathbf{Exp}(\widetilde{\theta })}$ satisfying conditions
\eqref{ast} and \eqref{martPPb}, and such that $Q$ is a PEMM for the process
$V$ with $V_{t}:=S_{t}-t\cdot \frac{\xi }{2}\cdot
\frac{\mathbb{E}_{P}[X_{1}\cdot e^{c\cdot X_{1}}]}{\mathbb{E}_{P}[e^{c\cdot
X_{1}}]} \in \mathcal{L}^{2}(Q)$ for any $t\geq 0$. But then, according to
Lemma~\ref{lem2} (a), we have
\begin{equation*}
Q_{X_{1}}(B)=\int_{B} \frac{(b-c)^{a}}{\Gamma (a)}
\cdot x^{a-1} \cdot e^{-(b-c)\cdot x}\,\lambda (dx)
\quad\text{for any } B\in \mathfrak{B}(\varUpsilon ).
\end{equation*}
The latter yields $Q_{X_{1}}=\mathbf{Ga}(\widetilde{\eta })$, where
$\widetilde{\eta }:=(b-c,a)\in \varUpsilon ^{2}$, and\enlargethispage{3pt}
\begin{equation*}
\mathbb{E}_{Q}[X_{1}]= \frac{\mathbb{E}_{P}[X_{1}\cdot e^{c\cdot
X_{1}}]}{\mathbb{E}_{P}[e^{c\cdot X_{1}}]}=
\frac{a}{b-c}>\frac{a}{b}=\mathbb{E}_{P}[X_{1}];
\end{equation*}
hence $p(P)<p(Q)<\infty $. In particular, according to Theorem~\ref{nfl} for
any $T>0$ the process $V_{\mathbb{T}}$ satisfies the (NFLVR) condition.
\end{ex}

%\begin{appendix}
%\end{appendix}

\begin{acknowledgement}[title={Acknowledgments}]
The authors would like to thank the anonymous reviewers for their valuable
comments and suggestions. Due to them, the presentation of the results is
more readable now.
\end{acknowledgement}\goodbreak

%\begin{funding}
%\gsponsor[id=,sponsor-id=]{}
%\gnumber[refid=]{}
%\end{funding}

\end{document}